\documentclass[12pt]{article}
\usepackage{amsmath,amssymb}

\newcommand{\Var}{{\cal{V}_{\mathbb{C}}}}

\def\tt{{\underline{T}}}

\def\nn{\underline{n}}

\def\kk{{\underline{k}}}

\def\1{\underline{1}}

\def\AA{{\mathbb A}}

\def\Z{{\mathbb Z}}

\def\H{{\mathbb H}}

\def\S{{\mathcal S}}

\def\ZZ{{\mathcal Z}}

\newtheorem{theorem}{Theorem}

\newenvironment{proof}
{\noindent{\bf Proof\/}.}{{ $\square$}\smallskip\par}

\title{Generating series of classes of Hilbert schemes of points on orbifolds
}
\author{S.M.~Gusein-Zade \thanks{Partially supported by the grants
RFBR-007-00593, INTAS-05-7805 and NWO-RFBR 047.011.2004.026.
Address: Moscow State University, Faculty
of Mathematics and Mechanics, Moscow, 119991, Russia. E-mail:
sabir\symbol{'100}mccme.ru} \and I.~Luengo \and
A.~Melle--Hern\'andez \thanks{The last two authors were partially
supported by the grant MTM2007-67908-C02-02. Address: University
Complutense de Madrid, Dept. of Algebra, Madrid, 28040, Spain.
E-mail: iluengo\symbol{'100}mat.ucm.es,
amelle\symbol{'100}mat.ucm.es}}

\date{}
\begin{document}
\def\eps{\varepsilon}

\maketitle

Let $K_0(\Var)$ be the Grothendieck
ring  of complex quasi-projective
varieties. This  is the abelian group generated by isomorphism classes $[X]$
of such varieties modulo the relation $[X]=[X-Y]+[Y]$ for a Zariski
closed subvariety $Y\subset X$; the multiplication is defined by the
Cartesian product: $[X_1]\cdot [X_2]=[X_1\times X_2]$.

Let $\mbox{Hilb}^n_X$, $n\ge 1$, be the Hilbert scheme of zero-dimensional subschemes
of length $n$ of a
complex
quasi-projective variety $X$. According to \cite{grot}, $\mbox{Hilb}^n_X$ is a quasi-projective variety. For a point $x\in X$,
let $\mbox{Hilb}^n_{X,x}$ be the Hilbert scheme of subschemes of $X$ supported
at the point $x$.

Let
$$
\H_X(T):=1+\sum\limits_{n=1}^\infty\,[\mbox{Hilb}^n_X]\,T^n \in 1+TK_0(\Var)[[T]],\quad \mbox{and}
$$
$$
\H_{X,x}(T):=1+\sum\limits_{n=1}^\infty\,[\mbox{Hilb}^n_{X,x}]\,T^n\in 1+TK_0(\Var)[[T]]\,\qquad
$$
be the generating series of classes of Hilbert schemes $\mbox{Hilb}^n_X$ and
$\mbox{Hilb}^n_{X,x}$ in the Grothendieck ring $K_0(\Var)$. 

\medskip

In \cite{MRL}, there was defined a notion of a power structure over a ring.
A \emph{power structure} over a commutative ring $R$ is a method to give sense to expressions of the form $(1+a_1T+a_2T^2+\ldots)^m$,
where $a_i$ and $m$ are elements of the ring $R$. In other words, a power structure is defined by a map
$\left(1+T\cdot R[[T]]\right)\times R\to 1+T\cdot R[[T]]$:
$$(A(T), m)\mapsto \left(A(T)\right)^m,\, (A(T)=1+a_1T+a_2T^2+\ldots,\,\,a_i\in R, m\in R),$$
such that all usual properties of the exponential
function hold.
In \cite{MRL}, there was described a natural power structure over
the Grothendieck ring $K_0(\Var)$ of complex quasi-projective varieties. 
It is closely connected with the $\lambda$-structure (see e.g. \cite{Knu}) on the Grothendieck ring $K_0(\Var)$ of quasi-projective varieties defined the Kapranov zeta function (\cite{Kap}).
The geometric description of this power structure is given as follows: if $A_1,A_2,\ldots,M$ are quasi-projective varieties, then the coefficient at $T^n$
in the series 
$$
(1+[A_1]T+[A_2]T^2+\ldots)^{[M]}
$$
is represented by the configuration space of pairs $(K,\varphi)$ consisting of  a finite subset $K$ of the variety $M$ and  a map $\varphi$  from $K$ to the disjoint
union $\coprod_{i=1}^\infty A_i$ of the varieties $A_i$, such that $\sum_{x\in K}I(\varphi(x))=n$, where $I:\coprod_{i=1}^\infty A_i\to {\Z}$ the tautological function sending the component $A_i$ to $i$.

There are two natural homomorphisms from the Grothendieck ring  $K_0(\Var)$ to the ring
$\Z$ of integers and to the ring $\Z[u,v]$ of polynomials in two variables: the Euler characteristic (with compact support) $\chi:K_0(\Var)\to \Z$ and the Hodge-Deligne polynomial
$e:K_0(\Var)\to \Z[u,v]$:
 $e(X)(u,v)=
\sum e^{p,q}(X) u^p v^q$. These homomorphisms respect the power structures on the correspondings rings
(see  e.g. \cite{nested}).

\medskip

For a smooth quasi-projective variety $X$ of dimension $d$, the following equation holds in $K_0(\Var)[[T]]$:
\begin{equation}\label{eq1}
\H_X(T)=\left(\H_{\AA^d,0}(T)\right)^{[X]}
\end{equation}
(\cite{Michigan}), where $\AA^d$
is the complex affine space of dimension $d$. For $d=2$ in other terms this was proved 
in the Grothendieck ring of motives by L.~G\"ottsche~\cite{Got2}.
For arbitrary $d$, the reduction of the equation (\ref{eq1}) under the Hodge-Deligne polynomial homomorphism $e:K_0(\Var)\to \Z[u,v]$ was proved by J.Cheah in \cite{Cheah}.

Since, for a point $x$ of a smooth variety $X$ of dimension $d$, the Hilbert scheme 
$\mbox{Hilb}^n_{X,x}$ can be identified with the Hilbert scheme
$\mbox{Hilb}^n_{\AA^d,0}$, the equation (\ref{eq1}) may be written as an integral with respect to the universal Euler characteristic $\chi_g(Y)=[Y]\in K_0(\Var)$ as follows

\begin{equation}\label{eq2}
\H_X(T):=\int_X \H_{X,x}(T)^{d\chi_g}.
\end{equation}
Here $d\chi_g$ is put in the exponent since the group operation in $1+T\cdot K_0(\Var)[[T]]$ is the multiplication.
Let $\psi$ be a constructible function on a quasi-projective variety $X$ with values in the abelian group $1+T\cdot K_0(\Var)[[T]]$. If $\psi$ is constant and equal to $\psi_\Sigma(T)$ on a stratum $\Sigma$ from a stratification $\S=\{\Sigma \}$ of $X$, then

\begin{equation}\label{eq3}
\int_X \H_{X,x}(T)^{d\chi_g}=\prod_{\Sigma\in \S} \left(\psi_{\Sigma}(T)\right)^{[\Sigma]}.
\end{equation}

If $X$ is not smooth but has isolated singularities, one can easlily see that the equation (\ref{eq2}) holds. However this colud not be the case in general.
Even the function $\H_{X,x}(T)$ (with values in $1+T\cdot K_0(\Var)[[T]]$)
could be not constructible: singularities of the space $X$ at points of some strata may have moduli.
It is interesting to understand to which extend the equation (\ref{eq2}) holds for varieties $X$ such that the function $\H_{X,x}(T)$ is constructible.

The described problem does not take place for orbifolds: they have finitely many local models. Here we prove the equation (\ref{eq2}) for complex (algebraic) orbifolds.

A $d$-dimensional \emph{orbifold} $X$ is a complex quasi-projective variety with an atlas of uniformizing systems for Zariski open sets in $X$ (see e.g. \cite{mp})). For a Zariski open subset
$U\subset X$, an \emph{uniformizing system} is a triple $({\widetilde U},G,\varphi)$, where $G$ is a finite group (depending on $U$), ${\widetilde U}$ is a smooth complex $d$-dimensional variety with a $G$-action , and $\varphi$ is an isomorphism ${\widetilde U}/G\to U$ of varieties (considered with the reduced structures).

For a point $x\in X$, let $({\widetilde U},G,\varphi)$ be an uniformizing system for a Zariski open neighbourhood  $U$ of the point $x$. Let $\pi_{{\widetilde U}}$ be the natural map ${\widetilde U} \to {\widetilde U}/G$ and let $\widetilde x\in (\varphi\circ\pi_{{\widetilde U}})^{-1}(x)$ be a representative of the corresponding orbit. The isotropy group $G_{\widetilde x}=\{g\in G: g\widetilde x=\widetilde x\,\}$ of the point $\widetilde x$ acts on the tangent space $T_{\widetilde x}{\widetilde U}$ by a representation $\alpha=\alpha_{\widetilde x}:G_{\widetilde x}\to GL(d,{\mathbb C})$. Moreover, there exists a system of local parameters $z_1,\ldots,z_d$ at the point $\widetilde x$ such that the action of the group $G_{\widetilde x}$ on the manifold ${\widetilde U}$ is given by standard linear equations corresponding to the representation:

\begin{equation}\label{eq4}
g^*z_i=\sum_{j=1}^d \alpha_{i,j}(g) z_j
\end{equation}
where $(\alpha_{i,j}(g))=\alpha(g)$.

One has: 1) the Hilbert scheme $\mbox{Hilb}^n_{X,x}$ of $0$-dimensional subschemes on $X$ supported at the point $x$ is isomorphic to the Hilbert scheme $\mbox{Hilb}^n_{\AA^d/{G_{\widetilde x}},0}$;
2) the partitioning of $U$ into parts corresponding to different conjugacy classes of isotropy subgroups $G_{\widetilde x}\subset G$ and to different (non-isomorphic) representations of the group $G_{\widetilde x}$ is a stratification of $U$. Therefore $\H_{X,x}(T)$ is a constructible function on $X$ with values in $1+T\cdot K_0(\Var)[[T]]$.

\begin{theorem}\label{theo1}
For an orbifold $X$, the following equation holds in $1+TK_0(\Var)[[T]]$:
\begin{equation}
\H_X(T)=\int_X \H_{X,x}(T)^{d\chi_g}.
\end{equation}
\end{theorem}

\begin{proof}
For a locally closed subvariety $Y\subset X$, let $\mbox{Hilb}^n_{X,Y}$
be the Hilbert scheme of subschemes of length $n$ of $X$ supported
at points of $Y$ and let
$$\H_{X,Y}(T):= 1+\sum\limits_{n=1}^\infty\,[\mbox{Hilb}^n_{X,Y}]\,T^k$$
be the corresponding generating series. If $Y$ is a Zariski closed subset
of $X$, then 

\begin{equation}\label{eq5}
\H_X(T)= \H_{X,Y}(T) \cdot \H_{X,X\setminus Y}(T).
\end{equation}

Therefore it is sufficient to prove the equation (\ref{eq2})
when $X$ is covered by one uniformizing system $({\widetilde U},G,\varphi)$ such that ${\widetilde U}$ lies in an affine space $\AA^N$.

Let us fix a subgroup $G'$ of the group $G$ and a representation $\alpha:G'\to GL(d,{\mathbb C})$. Let $X_{G',\alpha}$ be the image under $\varphi\circ \pi_G$ of the set of points $\widetilde x\in {\widetilde U}$ such that $G_{\widetilde x}=G'$ and the representation of the group $G'$ on the tangent space $T_{\widetilde x}{\widetilde U}$ is isomorphic to $\alpha$. Using (\ref{eq5}) again, one can see that it is sufficient to prove the formula 
\begin{equation}\label{eq6}
\H_{X,Y}(T)=\left(\H_{\AA^d/(G',\alpha),0}(T)\right)^{[Y]}
\end{equation}
for any point $x\in X_{G',\alpha}$ and for a certain Zariski open neighbourhood 
$Y$ of this point in an irreducible component of $X_{G',\alpha}$.

Let $\widetilde x\in (\varphi\circ\pi_{{\widetilde U}})^{-1}(x)$ be a point in ${\widetilde U}$
such that $G_{\widetilde x}=G'$. Then the representation of $G'$ on the tangent space $T_{\widetilde x}{\widetilde U}$ is automatically isomorphic to $\alpha$. Let $u_1,\ldots,u_d$ be a regular system of parameters at the point $\widetilde x$
(for example one may suppose that $\widetilde x=0$ and $u_1,\ldots,u_d$ are $d$ of the standard coordinates on $\AA^N$ such that the projection of the tangent space $T_{\widetilde x}{\widetilde U}$ to the corresponding $d$-dimensional coordinate plane is non-degenerate (i.e. $u_1-u_1^0,\ldots,u_d-u_d^0$ is a regular system of parameters at each point from a Zariski open neighbourhood of $\widetilde x$).

Let us suppose that the parameters  $u_1,\ldots,u_d$ are chosen in such a way that, in the corresponding coordinates on the tangent space $T_{\widetilde x}{\widetilde U}$, the action of the group $G'$ is given by the standard equations (\ref{eq4}). Define a new regular system of parameters $\widetilde u_1,\ldots,\widetilde u_d$ at the point $\widetilde x$ by the equations
$$
\widetilde u_i=\frac{1}{|G'|}\sum_{g\in G'}\sum_{j} \alpha_{i,j}(g^{-1})g^*u_j.
$$
One has 
$$
g^*\widetilde u_i=\sum_{j} \alpha_{i,j}(g)\widetilde u_j.
$$
Therefore, $\widetilde u_1-\widetilde u_1^0,\ldots,\widetilde u_d-\widetilde u_d^0$ is a regular system of parameters at each point from a Zariski open neighbourhood of $\widetilde x$ in the corresponding irreducible component of $(\varphi\circ \pi_G)^{-1} X_{G',\alpha}$. It defines a map $\widetilde{\underline u}-\widetilde{\underline u}^{\,0}:{\widetilde U},{\widetilde x'} \to \AA^d,0$.  There is a commutative diagram
\[\begin{array}{ccccc}  & {\widetilde U},{\widetilde x'} & \stackrel{\widetilde{\underline u}-\widetilde{\underline u}^{\,0}}{\longrightarrow} & \AA^d,0 & \\ &   \, \downarrow &&  \downarrow \,  &\\ & X,x' & \stackrel{}{\longrightarrow} & \AA^d/(G',\alpha), 0 \end{array}\]
and the lower map identifies the Hilbert scheme of points on $X$ supported at the point $x'$ with the Hilbert scheme of points on $\AA^d/(G',\alpha)$ supported at the origin.

This way, a zero-dimensional subscheme on $X$ supported at poins of $Y$ is defined by a finite subset $K\subset Y$ to each point $x$ of which there corresponds a zero-dimensional subscheme on $\AA^d/(G',\alpha)$ supported at the origin. The length of the subscheme is equal to the sum of lengths
of the corresponding subschemes of $\AA^d/(G',\alpha)$. As it follows from the geometric description of the power structure over the Grothendieck ring of quasi-projective varieties, the coefficient at $T^n$ in the right hand side of the equation (\ref{eq6}) is represented just by the configuration space of such objects. This proves the statement.
\end{proof}

Reductions of the equation (\ref{eq2}) under the Euler characteristic homomorphism and the 
Hodge--Deligne polynomial homomorphism gives equations 
for the generating series of the corresponding invariants of Hilbert schemes  of points on orbifolds. In particular
\begin{equation}
1+\sum_{n=1}^\infty\,\chi(\mbox{Hilb}^n_{X})T^n\,=\int_X (1+\sum_{n=1}^\infty\,\chi(\mbox{Hilb}^n_{X,x})T^n)^{d\chi}\,.
\end{equation}

\medskip

In \cite{Cheah2}, J.~Cheah considered nested Hilbert schemes on a smooth $d$-di\-men\-sio\-nal complex quasi-projective variety $X$.  For $\nn=(n_1,\ldots,n_r)\in\Z_{\geq 0}^r$, the \emph{nested Hilbert scheme} $Z^{\,\nn}_X$ of depth $r$ is the scheme which parametrizes collections of the form
$(Z_1,\ldots,Z_r)$, where $Z_i\in \mbox{Hilb}_X^{n_i}$ and $Z_i$ is a subscheme of $Z_j$ for $i<j$.
The scheme $Z^{\,\nn}_{X}$ is non-empty only if $n_1\leq n_2\leq \ldots\leq n_r$; notice that
$Z^{(n)}_{X}=\mbox{Hilb}_X^{n}\cong Z^{(n,\ldots, n)}_{X}$.

For $Y\subset X$,
let $Z^{\,\nn}_{X,Y}$ be the  scheme which parametrizes collections
$(Z_1,\ldots,Z_r)$ from $Z^{\nn}_X$ with $\mbox{supp}\,Z_i\subset Y$. For $Y=\{x\}$, $x\in X$, we shall use the notation $Z^{\,\nn}_{X,x}$.

For $r\geq 1$, let $\tt=(T_1,\ldots,T_r)$ and
$$
\ZZ_X^{(r)}(\tt):=1+\sum\limits_{\nn\in \Z_{\geq 0}^r\setminus\{ 0\}}\,[Z^{\,\nn}_X]\,\tt^{\nn},\quad
\ZZ_{X,x}^{(r)}(\tt):=1+\sum\limits_{\nn\in \Z_{\geq 0}^r\setminus\{ 0\}}\,[Z^{\,\nn}_{X,x}]\,\tt^{\nn},\quad
$$
be the generating series of classes of the nested Hilbert schemes $Z^{\nn}_X$ of depth $r$ (resp. of those supported at the point $x$).

A series $A(\tt)=1+\sum\limits_{\nn\in \Z_{\geq 0}^r\setminus\{ 0\}}\,a_{\nn} \,\tt^{\nn},$
${\tt}^{\nn}:=T_1^{n_1}\cdot\ldots\cdot T_r^{n_r}$, $a_{\nn}\in K_0(\Var)$, has a unique representation of the form 
$A(\tt)=\prod\limits_{\kk\in \Z_{\geq 0}^r\setminus\{ 0\}} (1-\tt^{\kk})^{-s_\kk}$, 
$s_\kk\in K_0(\Var)$. Then for $m\in K_0(\Var)$, $A(\tt)^m:=\prod\limits_{\kk\in \Z_{\geq 0}^r\setminus\{ 0\}} (1-\tt^{\kk})^{-ms_\kk}$ (see \cite{nested}).

For a smooth quasi-projective variety $X$ of dimension $d$, the following equation holds in  $K_0(\Var)[[\tt]]${}$)$:
\begin{equation*}
\ZZ_X^{(r)}(\tt)=\left(\ZZ_{\AA^d,0}^{(r)}(\tt)\right)^{[X]}
\end{equation*}
(\cite{nested}).

Using the same arguments as in the proof of Theorem 1, one gets the following statetment

\begin{theorem}\label{theo2}
For an  orbifold $X$, the following equation holds
\begin{equation*}\label{eq7}
\ZZ_X^{(r)}(\tt)=\int_X \left(\ZZ_{X,x}^{(r)}(\tt)\right)^{d\chi_g}.
\end{equation*}
\end{theorem}

\end{document}